  \documentclass{monsky2009} 

\usepackage{graphicx}

\usepackage{amssymb}
\usepackage{amsmath}

\usepackage{bm}

\newenvironment{theorem*}[1]{\textbf{#1}\itshape \hspace{.3em}}{\upshape}
\newenvironment{remark*}[1]{\textbf{#1}\itshape \hspace{.3em}}{\upshape}
\newenvironment{corollary*}[1]{\textbf{#1}\itshape \hspace{.3em}}{\upshape}
\newenvironment{proof}{\textbf{Proof\hspace{.3em}}}{}

\newtheorem{definition}{Definition}[section]
\newtheorem{theorem}[definition]{Theorem}
\newtheorem{lemma}[definition]{Lemma}

\newtheorem{corollary}[definition]{Corollary}

\newcommand{\charstic}{\ensuremath{\mathrm{char\ }}}

\newcommand{\dif}{\ensuremath{\mathrm{d}}}

\newcommand{\bilprod}{\mbox{\ensuremath\,\#\,}}

\begin{document}
\addtolength{\topmargin}{-30pt}
\addtolength{\textheight}{30pt}

\begin{frontmatter}






\title{Transcendence of some Hilbert-Kunz multiplicities (modulo a conjecture)}
\author{Paul Monsky}

\address{Brandeis University, Waltham MA  02454-9110, USA. monsky@brandeis.edu}

\begin{abstract}
Suppose that $h\in F[x,y,z]$, $\charstic F=2$, defines a nodal cubic. In earlier papers we made a precise conjecture as to the Hilbert-Kunz functions attached to the powers of $h$. Assuming this conjecture we showed that a class of characteristic 2 hypersurfaces has algebraic but not necessarily rational Hilbert-Kunz multiplicities. We now show that if the conjecture holds, then transcendental multiplicities exist, and in particular that the number $\sum \binom{2n}{n}^{2}/(65,536)^{n}$, proved transcendental by Schneider, is a $Q$-linear combination of Hilbert-Kunz multiplicities of characteristic 2 hypersurfaces.
\end{abstract}


\end{frontmatter}


\section{The power series $\theta_{g}$}
\label{section1}

We recall some notation and results from \cite{3} that will be used throughout. $X$ is the vector space of functions $I\rightarrow Q$ where $I=[0,1]\cap Z[\frac{1}{2}]$. $F$ is a field of characteristic 2, while $q$ is always a power, $2^{n}$, of 2. For each $f\ne 0$ in the maximal ideal of $F[[x_{1},\ldots , x_{r}]]$ there is an element $\phi_{f}$ of $X$ whose value at $\frac{i}{q}$ is $q^{-r}e_{n}(f^{i})$ where $n\rightarrow e_{n}$ is the Hilbert-Kunz function. There is a symmetric bilinear product $\bilprod$, $X\times X \rightarrow X$ with the following property. If $f$ and $g$ are in $F[[u_{1},\ldots , u_{r}]]$ and $F[[v_{1},\ldots , v_{s}]]$, and $h=f+g$, then $\phi_{h}=\phi_{f}\bilprod \phi_{g}$. $\Gamma$ is the Grothendieck group of isomorphism classes of finitely generated $F[T]$-modules annihilated by a power of $T$; see the material following Theorem 1.5 of \cite{3} for the $Z$-basis $\lambda_{0},\lambda_{1},\ldots$ of $\Gamma$ and the multiplication on $\Gamma$. For $n \ge 0$ and $\alpha$ in $X$, $\mathcal{L}_{n}(\alpha)$ is the element $\sum_{0}^{q-1}\left(\alpha\left(\frac{i+1}{q}\right)-\alpha\left(\frac{i}{q}\right)\right)(-)^{i}\lambda_{i}$ of $\Gamma_{Q}=\Gamma\underset{Z}{\otimes}Q$. The product in $\Gamma_{Q}$ of $\mathcal{L}_{n}(\alpha)$ and $\mathcal{L}_{n}(\beta)$ is $\mathcal{L}_{n}(\alpha\bilprod\beta)$.

\begin{definition}
\label{def1.1}
$\eta: \Gamma_{Q}\rightarrow Q$ is the linear map taking each $\lambda_{i}$ to $1$.
\end{definition}

\begin{lemma}
\label{lemma1.2}
If $\gamma= \alpha\bilprod\beta$ then for each $n$, $\eta(\mathcal{L}_{n}(\gamma))=\eta(\mathcal{L}_{n}(\alpha))\cdot\eta(\mathcal{L}_{n}(\beta))$.
\end{lemma}

\begin{proof}
If $k$ is the Nim-sum of $i$ and $j$ then $\lambda_{i}\cdot\lambda_{j}=\lambda_{k}$. It follows that $\eta(\lambda_{i}\lambda_{j})$ and $\eta(\lambda_{i})\eta(\lambda_{j})$ are both $1$ and that $\eta: \Gamma_{Q}\rightarrow Q$ is multiplicative. Now use the fact that $\mathcal{L}_{n}(\gamma)=\mathcal{L}_{n}(\alpha)\mathcal{L}_{n}(\beta)$.
\qed
\end{proof}

\begin{lemma}
\label{lemma1.3}
If $n\ge 1$, $\eta(\mathcal{L}_{n}(\alpha))=2^{2n}(\phi_{uv}\bilprod \alpha)\left(\frac{1}{2^{n}}\right)-2^{2n-1}(\phi_{uv}\bilprod \alpha)\left(\frac{1}{2^{n-1}}\right)$.
\end{lemma}

\begin{proof}
Write $\mathcal{L}_{n}(\alpha)$ as $\sum_{0}^{q-1}(-)^{i}a_{i}\lambda_{i}$ where $q=2^{n}$. Now $2^{2n}(\phi_{uv}\bilprod \alpha)\left(\frac{1}{2^{n}}\right)$ is the co-efficient of $\lambda_{0}$ in $2^{2n}\mathcal{L}_{n}(\phi_{uv}\bilprod \alpha)=2^{2n}\mathcal{L}_{n}(\phi_{uv})\mathcal{L}_{n}(\alpha)$. Since $e_{n}(u^{i}v^{i})=q^{2}-(q-i)^{2}$, $2^{2n}\mathcal{L}_{n}(\phi_{uv})=\sum_{0}^{q-1}(-)^{i}(2q-1-2i)\lambda_{i}$. Using the formula for the product in $\Gamma$ we find that $2^{2n}(\phi_{uv}\bilprod \alpha)\left(\frac{1}{2^{n}}\right)$ is $\sum_{0}^{q-1}(2q-1-2i)a_{i}$.

Rewrite the above as $\sum_{0}^{\frac{q}{2}-1}(2q-1-4i)(a_{2i}+a_{2i+1})-2(a_{1}+a_{3}+\cdots +a_{q-1})$. Now $\mathcal{L}_{n-1}(\alpha)$ is easily seen to be $\sum_{0}^{\frac{q}{2}-1}(-)^{i}(a_{2i}+a_{2i+1})\lambda_{i}$. The last paragraph, with $n$ replaced by $n-1$, shows that $2^{2n-2}(\phi_{uv}\bilprod \alpha)\left(\frac{1}{2^{n-1}}\right)=\sum_{0}^{\frac{q}{2}-1}(q-1-2i)(a_{2i}+a_{2i+1})$. It follows that the right-hand side of Lemma \ref{lemma1.3} is $\sum_{0}^{q-1}a_{i}-2(a_{1}+a_{3}+\cdots + a_{q-1})=\eta(\mathcal{L}_{n}(\alpha))$.
\qed
\end{proof}

\begin{corollary}
\label{corollary1.4}
Let $g$ be a power series in $r$ variables over $F$, and $G=uv+g$, viewed as a power series in $r+2$ variables. Then for $n\ge 1$, $e_{n}(G)-2^{r+1}e_{n-1}(G)=2^{rn}\eta(\mathcal{L}_{n}(\phi_{g}))$.
\end{corollary}

\begin{proof}
Set $\alpha=\phi_{g}$, and multiply the identity of Lemma \ref{lemma1.3} by $2^{rn}$. Then since $\phi_{uv}\bilprod \alpha = \phi_{uv}\bilprod\phi_{g} = \phi_{G}$, the right-hand side becomes $e_{n}(G)-\frac{1}{2}\cdot 2^{r+2}\cdot e_{n-1}(G)$, giving the corollary.
\qed
\end{proof}

\begin{definition}
\label{def1.5}
Let $r$, $g$ and $G$ be as in Corollary \ref{corollary1.4}.  Then $\theta_{g}$ in $Z[[w]]$ is $(1-2^{r+1}w)\sum{e_{n}}(G)w^{n}$.
\end{definition}

Note that the power series $\theta_{g}$ converges in $|w|\le \frac{1}{2^{r}}$, and that $\theta_{g}\left(\frac{1}{2^{r+1}}\right)$ is just the Hilbert-Kunz multiplicity of $G$. Furthermore the co-efficient of $w^{n}$ in $\theta_{g}$ is $e_{n}(G)-2^{r+1}e_{n-1}(G)$.

Suppose now that we have finitely many power series $g_{i}$ over $F$, that $g_{i}$ is in $r_{i}$ variables, and that the variables corresponding to distinct $g_{i}$ are pairwise disjoint.

\begin{theorem}
\label{theorem1.6}
$\theta_{\Sigma g_{i}}$  is the Hadamard product of the $\theta_{g_{i}}$.
\end{theorem}

\begin{proof}
Let $G_{i}=uv+g_{i}$ and $M=uv+\sum g_{i}$. Applying Corollary \ref{corollary1.4} to each $g_{i}$ and to $\sum g_{i}$ and using Lemma \ref{lemma1.2} repeatedly we find that $e_{n}(M)-2^{1+\Sigma r_{i}}e_{n-1}(M)$ is the product of the various $e_{n}(G_{i})-2^{1+r_{i}}e_{n-1}(G_{i})$. In other words the co-efficient of $w^{n}$ in $\theta_{\Sigma g_{i}}$ is the product of the co-efficients of $w^{n}$ in the various $\theta_{g_{i}}$. 
\qed
\end{proof}

\begin{corollary}
\label{corollary1.7}
The Hilbert-Kunz multiplicity of $uv+\sum g_{i}$ is the value at $w=\frac{1}{2^{1+\Sigma r_{i}}}$ of the Hadamard product of the $\theta_{g_{i}}$.
\end{corollary}

Corollary \ref{corollary1.7} is the key to this note. It allows us to pass from the algebraic realm (the $\theta$ attached to $x^{3}+y^{3}+xyz$ is conjectured by us to be algebraic of degree 2 over $Q(w)$) to the transcendental realm by making use of the Hadamard product. We next give some easy results about $\theta$.

\begin{theorem}
\label{theorem1.8}
Let $g$ be an $r$-variable power series and $h=g^{2}$. Then $\theta_{h}=1+2^{r}w\cdot \theta_{g}$.
\end{theorem}

\begin{proof}
If $i\le q$, $\phi_{h}\left(\frac{i}{2q}\right)=\phi_{g}\left(\frac{2i}{2q}\right)$, while if $i\ge q$, $\phi_{h}\left(\frac{i}{2q}\right)=1$. It follows that $\mathcal{L}_{n+1}(\phi_{h})=\mathcal{L}_{n}(\phi_{g})$ for all $n$. Applying $\eta$ and multiplying by $2^{r(n+1)}$ we find from Corollary \ref{corollary1.4} that the left hand side of this equation becomes $e_{n+1}(H)-2^{r+1}e_{n}(H)$ where $H=uv+h$. And the right-hand side becomes $2^{r}(e_{n}(G)-2^{r+1}e_{n-1}(G))$ where $G=uv+g$. We have shown that the co-efficients of $w^{n+1}$ in $\theta_{h}$ and in $2^{r}w\cdot \theta_{g}$ are equal, giving the theorem.
\qed
\end{proof}

When $r=1$ and $g$ is a power of the variable it's easy to calculate $\theta_{g}$. In particular we find:

\begin{lemma}
\label{lemma1.9}
When $r=1$,
\begin{enumerate}
\item[(a)] $\theta_{x^{5}}-\theta_{x^{3}} = 2\sum_{1}^{\infty}w^{2n}$
\item[(b)] $2\theta_{x^{3}}-\theta_{x^{5}} = 1+2\sum_{0}^{\infty}w^{2n+1}$
\end{enumerate}
Now take $r=3$, fix $f$ in $F[[x,y,z]]$, and let $\theta_{f}=\sum a_{n}w^{n}$. We shall use Corollary \ref{corollary1.7} to show that certain infinite sums involving the $a_{n}$ are $Q$-linear combinations of Hilbert-Kunz multiplicities.
\end{lemma}

\begin{theorem}
\label{theorem1.10}
$\sum\frac{a_{2n}^{2}}{2^{16n}}$, $\sum\frac{a_{2n+1}^{2}}{2^{16n}}$, $\sum\frac{a_{2n}a_{2n+1}}{2^{16n}}$, $\sum\frac{a_{2n+1}a_{2n+2}}{2^{16n}}$ and $\sum\frac{a_{2n}a_{2n+2}}{2^{16n}}$ are all $Q$-linear combinations of Hilbert-Kunz multiplicities.
\end{theorem}

\begin{proof}
Using Theorem \ref{theorem1.8} we find that the co-efficients of $w^{2n+2}$ in $\theta_{f}$ and $\theta_{f^{2}}$ are $a_{2n+2}$ and $8a_{2n+1}$. So the Hadamard product of $\theta_{f}$, $\theta_{f^{2}}$ and $\theta_{x^{5}}-\theta_{x^{3}}=2\sum_{0}^{\infty}w^{2n+2}$ is $16\sum_{0}^{\infty}a_{2n+1}a_{2n+2}w^{2n+2}$. Since $1+(3+3+1)=8$, Corollary \ref{corollary1.7} with $g_{1}=f$, $g_{2}=f(X,Y,Z)^{2}$ and $g_{3}=T^{3}$ or $T^{5}$ shows us that this Hadamard product, evaluated at $\frac{1}{2^{8}}$, is a difference of two Hilbert-Kunz multiplicities. Also, the co-efficients of $w^{2n+1}$ in $\theta_{f}$ and $\theta_{f^{2}}$ are $a_{2n+1}$ and $8a_{2n}$. So the Hadamard product of $\theta_{f}$, $\theta_{f^{2}}$ and $2\theta_{x^{3}}-\theta_{x^{5}}=1+2\sum w^{2n+1}$ is $1+16\sum_{0}^{\infty}a_{2n}a_{2n+1}w^{2n+1}$. As above we see that this Hadamard product, evaluated at $\frac{1}{2^{8}}$, is a $Z$-linear combination of two Hilbert-Kunz multiplicities. We have proved the third and fourth of the assertions of the theorem. To prove the first and the second we argue similarly with $g_{1}=f$ and $g_{2}=f(X,Y,Z)$. For the final assertion we take $g_{1}=f$ and $g_{2}=f(X,Y,Z)^{4}$, using Theorem \ref{theorem1.8} to calculate $\theta_{g_{2}}$.
\qed
\end{proof}

\section{Transcendence results (modulo a conjecture)}
\label{section2}

Now let $r=3$, and $f=x^{3}+y^{3}+xyz$ be the defining equation of a nodal cubic. In Definitions 2.1 and 2.2 of \cite{3} we constructed elements $\phi_{0},\phi_{1},\ldots$ of $X$; we further conjectured that $\phi_{f}=t+\phi_{0}$. (This is an alternative form of a conjecture we made earlier in \cite{2}.)

\begin{lemma}
\label{lemma2.1}
Let $A_{n}$ be the binomial co-efficient $\binom{2n}{n}$, so that $\sum A_{n}w^{2n}$ converges to $(1-4w^{2})^{-\frac{1}{2}}$ in the disc $|w|<\frac{1}{2}$. Then if the conjecture of \cite{3} holds, $(1-6w+8w^{2})\theta_{f}=(1-2w)+(2w-8w^{2}-24w^{3})\sum A_{n}w^{2n}$.
\end{lemma}

\begin{proof}
This is shown in the paragraph following Corollary 2.7 of \cite{2}. Alternatively it is an easy consequence of Lemma 2.6 of \cite{3}, since in the notation of that lemma, $\sum e_{n}(uv+f)w^{n}$ is, granting the conjecture, equal to $\sum (2^{-n}+E_{1}(2^{-n}))(32w)^{n}$.
\qed
\end{proof}

\begin{theorem}
\label{theorem2.2}
If the conjecture of \cite{3} holds, then the value of $\sum A_{n}^{2}\lambda^{2n}$ at $\lambda = \frac{1}{2^{8}}$ is a $Q$-linear combination of Hilbert-Kunz multiplicities.
\end{theorem}

\begin{proof}
Suppose $\theta_{f}=\sum a_{n}w^{n}$. Comparing co-efficients of $w^{2n+2}$ in the identity of Lemma \ref{lemma2.1} we see that $a_{2n+2}-6a_{2n+1}+8a_{2n}=-8A_{n}$. So $64\sum A_{n}^{2}\cdot\left(\frac{1}{2^{16n}}\right) = \sum (a_{2n+2}-6a_{2n+1}+8a_{2n})^{2}\cdot\left(\frac{1}{2^{16n}}\right)$. Expanding and using Theorem \ref{theorem1.10} we get the result.
\qed
\end{proof}

It only remains to show that $\sum A_{n}^{2}\cdot\left(\frac{1}{2^{16n}}\right)$ is transcendental. Results of this sort about special values of hypergeometric functions were first shown by Schneider \cite{4}; we'll sketch a proof.

\begin{lemma}
\label{lemma2.3}
\textbf{\upshape(Euler)} If $0\le \lambda < \frac{1}{4}$, $\int_{-1}^{1}\frac{\dif w}{\sqrt{(1-16\lambda^{2}w^{2})(1-w^{2})}} = \pi\sum A_{n}^{2}\lambda^{2n}$.
\end{lemma}

\begin{proof}
$(1-16\lambda^{2}w^{2})^{\frac{1}{2}}=\sum(4\lambda^{2})^{n}A_{n}w^{2n}$. So our integral is $\sum (4\lambda^{2})^{n}A_{n}\cdot\int_{-1}^{1}\frac{w^{2n}\dif w}{\sqrt{1-w^{2}}}$. But as one learns in every introductory calculus course, $\int_{-1}^{1}\frac{w^{2n}\dif w}{\sqrt{1-w^{2}}}=\frac{\pi A_{n}}{4^{n}}$.
\qed
\end{proof}

Now suppose that $\lambda$ in $(0,\frac{1}{4})$ is rational, and consider the affine curve $y^{2}=(1-16\lambda^{2}x^{2})(1-x^{2})$ defined over $Q$.  The real locus of this curve has a single bounded component. This component is a simple closed curve containing $(-1,0)$ and $(1,0)$. $\frac{\dif x}{y}$ is a $1$-form on our affine curve. When we integrate $\frac{\dif x}{y}$ over the bounded component (with the clockwise orientation) the integrals over the pieces in the upper and lower half-planes are equal, and each is $\int_{-1}^{1}\frac{\dif x}{\sqrt{(1-16\lambda^{2}x^{2})(1-x^{2})}}$. So by Lemma \ref{lemma2.3} the integral of $\frac{\dif x}{y}$ over the bounded component is $2\pi \sum A_{n}^{2}\lambda^{2n}$.

\begin{lemma}
\label{lemma2.4}
Suppose $\lambda$ in $(0,\frac{1}{4})$ is rational. Then there is an elliptic curve $E$ defined over $Q$, a holomorphic $1$-form on $E$ rational over $Q$ and a $1$-cycle on $E$, such that the integral of the form over the $1$-cycle is $2\pi \sum A_{n}^{2}\lambda^{2n}$.
\end{lemma}

\begin{proof}
Projectify the affine curve  $y^{2}=(1-16\lambda^{2}x^{2})(1-x^{2})$, and let $E$ be a desingularization over $Q$ of the resulting complex projective plane quartic. The $1$-form is the pull-back of $\frac{\dif x}{y}$, and the $1$-cycle corresponds to the inverse image in $E$ of the bounded component of the real locus of the affine curve $y^{2}=(1-16\lambda^{2}x^{2})(1-x^{2})$.
\qed
\end{proof}

\begin{theorem}
\label{theorem2.5}
\textbf{\upshape(Schneider)} Let $E$ be an elliptic curve defined over $Q$ and $\omega\ne 0$ be the integral over a $1$-cycle on $E$ of a holomorphic $Q$-rational $1$-form. The $\frac{\omega}{\pi}$ is transcendental.
\end{theorem}

\begin{proof}
On page 63 of \cite{4}, Schneider deduces this from his Theorem 13. Lang later put Theorem 13 into a more general setting. I'll indicate how Theorem \ref{theorem2.5} above follows from this ``Schneider-Lang theorem''.

We may assume that the curve is $y^{2}=4x^{3}-g_{2}x-g_{3}$ with $g_{2}$ and $g_{3}$ rational, and that the $1$-form is $\frac{\dif x}{y}$. Let $\wp$ be the associated Weierstrass $\wp$-function. $\wp(z+\omega)=\wp(z)$ and $(\wp^{\prime})^{2}=4\wp^{3}-g_{2}\wp-g_{3}$. Choose $z_{0}$ with $\wp(z_{0})$ and consequently $\wp^{\prime}(z_{0})$ algebraic. Suppose that $\frac{\omega}{\pi}$ is algebraic. We define a subring $Y$ of the field of meromorphic functions on $\mathbb{C}$ --- $Y$ is generated over $Q$ by $\wp(z+z_{0})$, $\wp^{\prime}(z+z_{0})$, $\mathrm{e}^{\frac{2\pi \mathrm{i} z}{\omega}}$, and the algebraic constants $\wp(z_{0})$, $\wp^{\prime}(z_{0})$ and $\frac{2\pi \mathrm{i}}{\omega}$. $Y$ is stable under $\frac{\dif}{\dif z}$. Furthermore the values of the 6 generators of $Y$ at the infinitely many points $\omega,2\omega,3\omega,\ldots$ all lie in a fixed number field. The Schneider-Lang theorem \cite{1} then tells us that the transcendence degree of $Y$ over $Q$ is $\le 1$ which is evidently false.
\qed
\end{proof}

Combining Lemma \ref{lemma2.4} with Theorem \ref{theorem2.5} we find:

\begin{corollary}
\label{corollary2.6}
If $\lambda$ in $(0,\frac{1}{4})$ is rational then $\sum A_{n}^{2}\lambda^{2n}$ is transcendental.
\end{corollary}

So if the conjecture of \cite{3} holds, the transcendental number $\sum \frac{A_{n}^{2}}{2^{16n}}$ is a $Q$-linear combination of Hilbert-Kunz multiplicities, and transcendental Hilbert-Kunz multiplicities exist.



\label{}



\end{document}